\documentclass[11pt]{article}
\usepackage{amsmath,amsfonts}
\usepackage{verbatim}
\usepackage{relsize}
\usepackage{enumerate}
\usepackage[inline]{enumitem}
\normalbaselines
\DeclareMathOperator{\diver}{div}

\newtheorem{theorem}{Theorem}
\newtheorem{proposition}[theorem]{Proposition}
\newtheorem{lemma}[theorem]{Lemma}
\newtheorem{corollary}[theorem]{Corollary}
\newtheorem{remark}[theorem]{Remark}

\newenvironment{proof}{\rm \trivlist \item[\hskip \labelsep{\it
      Proof}.]}{\nopagebreak \hfill $\square$ \endtrivlist}

\renewcommand{\div}{\diver}
\newcommand{\gbar}{\overline{g}}
\newcommand{\Mbar}{\overline{M}}
\newcommand{\Ubar}{\overline{U}}
\newcommand{\Xbar}{\overline{X}}
\newcommand{\nablabar}{\overline{\nabla}}
\newcommand{\divbar}{\overline{\div}}

\newcommand{\PP}{\mathbb{P}}
\newcommand{\R}{\mathbb{R}}

\newcommand{\IxP}{I {}_h\hspace{-2pt}\times_{\rho} \PP}
\newcommand{\di}{\mathrm{d}}

\newcommand{\trace}{\mathrm{trace}}
\newcommand{\II}{\mathrm{I\hspace{-1pt}I}}

\title{A note on spacelike hypersurfaces and timelike conformal vectors}
\author{Giulio Colombo${}^{a}$, Jos\'e A. S. Pelegr\'in${}^{b}$ and Marco Rigoli${}^{a}$ \\[6mm]
${}^a$ Dipartimento di Matematica, Universit\`a degli Studi di Milano, \\[0.5mm] via Cesare Saldini 50, 20133 Milano, Italy,\\[0.5mm] E-mail\textup{: \texttt{giulio.colombo@unimi.it, 
marco.rigoli@unimi.it}}\\[3mm]
${}^b$ Departamento de Did\'actica de la Matem\'atica, \\[0.5mm] Universidad de Granada, 18071 Granada, Spain \\ E-mail\textup{:
\texttt{jpelegrin@ugr.es}} \\[3mm]}

\date{}
\begin{document}

\maketitle

\thispagestyle{empty}

\abstract{Any compact spacelike hypersurface immersed in a doubly warped product spacetime $\IxP$ with nondecreasing warping factor $\rho$ must be a spacelike slice, provided that the mean curvature satisfies $H\geq\rho'/h\rho$ everywhere on the hypersurface. The conclusion also holds, under suitable assumptions on the immersion, when the hypersurface is complete and noncompact. A similar rigidity property is shown for compact hypersurfaces in spacetimes carrying a conformal, strictly expanding, timelike vector field.}

\section{Introduction}

Spacelike hypersurfaces play a crucial role in the understanding of the geometry of a Lorentzian spacetime. Roughly speaking, they describe the physical space that can be measured in a given instant of time. For instance, they serve as initial data in the Cauchy problem for Einstein's field equations \cite {Ri} and they play a privileged role in determining the causal properties of the spacetime. Indeed, a spacetime is globally hyperbolic if and only if it admits a Cauchy hypersurface \cite{Ger}. Even more, any globally hyperbolic spacetime is diffeomorphic to $\mathbb{R} \times S$, being $S$ a smooth spacelike Cauchy hypersurface \cite{BS}.

In this article we will study the geometry of spacelike hypersurfaces in certain spacetimes that present a particular symmetry. In General Relativity, symmetry arises from the existence of a one-parameter group of transformations generated by a conformal vector field. This infinitesimal symmetry is usually assumed when searching for exact solutions of Einstein's field equations \cite{Ea}. In fact, there is a vast literature concerning the study of spacelike hypersurfaces in spacetimes that admit different causal symmetries (see \cite{CRR}, \cite{CPR}, \cite{CPR2}, for instance). We will focus on the study of spacelike hypersurfaces in doubly warped product spacetimes, where natural conformal vector field pertains to this structure.

Let $(\PP,\sigma)$ be a (connected) Riemannian manifold of dimension $m$ and $I\subseteq\mathbb R$ be an open interval. The product manifold $\Mbar = I \times \PP$ can be endowed with the Lorentzian metric $\gbar$ given at each point $(t,x)\in\Mbar$ by 
\begin{equation}
	\gbar_{(t,x)} = -h(x)^2\pi_I^{\ast}(\di t^2) + \rho(t)^2\pi_{\PP}^{\ast}(\sigma),
\end{equation}
where $h\in C^{\infty}(\PP)$, $\rho\in C^{\infty}(I)$ are positive functions and $\pi_I:\Mbar\to I$, $\pi_{\PP}:\Mbar\to \PP$ are the projections onto the factors of $I\times \PP$. We call $(\Mbar, \gbar)$ a doubly warped product spacetime, which will be denoted here by $\Mbar = \IxP$. The time orientation of $\Mbar$ is the one given by the timelike vector field $\partial_t := \partial /\partial t$. The causal symmetry in these ambient spacetimes is given by the timelike vector field $\rho \partial_t$, which is conformal, as we will prove in Lemma \ref{lemmaconfo} below. As a consequence, the family of spacelike slices $\Sigma_t = \{t\}\times\PP$, $t\in I$, provides a foliation of $\Mbar$ by totally umbilical hypersurfaces; in particular, the mean curvature of $\Sigma_t$ at the point $(t,x) \in \Mbar$ is given by $\mathcal H(t,x) = \rho'(t)/\rho(t)h(x)$.

Doubly warped product spacetimes were introduced in \cite{BP} as an extension of the singly warped product spacetimes defined in \cite{O'N} and include several models such as standard static spacetimes and Generalized Robertson-Walker spacetimes (see \cite{BO} for the original definition of singly warped product spaces in the Riemannian setting). Causal properties of doubly warped product spacetimes were studied in \cite{Al}, obtaining several conditions to guarantee their global hyperbolicity. Moreover, in \cite{GO} the authors studied the decomposition of an ambient spacetime as a doubly warped product.

Our article is organized as follows. Section \ref{se1} is devoted to proving some preliminary results for spacelike hypersurfaces in general Lorentzian manifolds as well as for the particular case of doubly warped product spacetimes. These results will enable us to prove our main theorems in Section \ref{se2}. In particular, we prove in Theorem \ref{thcpt} that a compact spacelike hypersurface $\psi : M \to \Mbar$ immersed in a doubly warped product spacetime $\Mbar = \IxP$ must be a spacelike slice if $H \geq \mathcal H\circ\psi \geq 0$ on $M$, with $H$ the mean curvature of $\psi$ in the direction of the normal unit vector with the same time-orientation as $\partial_t$. Theorem \ref{thcptgen} partially extends this observation to the case of a general spacetime carrying a conformal timelike vector field. Theorems \ref{thasymp}, \ref{thparab1} and \ref{thparab2} provide examples of cases where the conclusion of Theorem \ref{thcpt} still holds true for $M$ complete and non-compact, under additional assumptions on $\psi$.

\section{The geometric setting}

\label{se1}

Let $\psi:(M^m,g)\to(\Mbar^{m+1},\gbar)$ be an isometric immersion between (connected) semi-Riemannian manifolds and denote with $\nabla$ and $\nablabar$, respectively, the Levi-Civita connections for $g$ and $\gbar$. For any given point $p\in M$, there exist sufficiently small neighbourhoods $U\subseteq M$ and $\Ubar\subseteq\Mbar$, respectively, of $p$ and $\psi(p)$, such that $\psi(U)\subseteq\Ubar$, $\psi|_U : U \to \Ubar$ is an embedding and $\Ubar$ supports a nowhere vanishing vector field $Z\in\mathfrak{X}(\Ubar)$ satisfying $\gbar(Z_{\psi(q)},(\di\psi)_q V)=0$ for every $q\in U$, $V\in T_q M$. We say that $Z$ is orthogonal to $\psi$ on $U$. Note that we are not requiring $\gbar(Z,Z)$ to be constant on $\Ubar$. From now on, for the sake of simplicity, we shall omit writing pullbacks or pushforwards via $\psi$ since the embedding $\psi|_U$ provides a natural identification between $U$ and a regular submanifold of $\Ubar$. Following Chapter 4 of \cite{O'N}, the second fundamental form $\II$ of $\psi$ is defined by
\begin{equation}
    \II(V,W) = \nablabar_V W - \nabla_V W
\end{equation}
for every couple of vectors $V,W$ tangent to $U$. $\II$ is a symmetric bilinear form taking values in the normal bundle $TU^{\bot} \subseteq T\Ubar$ of vectors orthogonal to $U$ and satisfies Weingarten's equation
\begin{equation} \label{wein}
    \gbar(\II(V,W),Z) = \gbar(\nablabar_V W, Z) = - \gbar(\nablabar_V Z, W)
\end{equation}
For $\Xbar\in\mathfrak X(M)$, $V,W\in T_q\Mbar$, $q\in\Mbar$, we have the differential identity
\begin{equation} \label{lieg}
    (\mathcal L_{\Xbar}\gbar)(V,W) = \gbar(\nablabar_V\Xbar,W) + \gbar(\nablabar_W\Xbar,V)
\end{equation}
for the Lie derivative $\mathcal L$ of the metric $\gbar$, so equation \eqref{wein} can be rewritten as
\begin{equation} \label{IIlie}
    \gbar(\II(V,W),Z) = - \frac{1}{2}(\mathcal L_Z\gbar)(V,W).
\end{equation}
The mean curvature vector $\mathbf{H}$ of $\psi$ is the normalized trace of $\II$ with respect to the metric $g$, that is, for every $q\in U$ and for every $g$-orthonormal basis $\{e_i\}_{1\leq i\leq m}$ of $T_q U$,
\begin{equation} \label{Hdef}
    \mathbf{H}_q = \frac{1}{m}\trace_g(\II_q) = \frac{1}{m}\sum_{i=1}^m g(e_i,e_i)\II(e_i,e_i)
\end{equation}

For a generic vector field $\Xbar\in\mathfrak X(\Mbar)$ we let $X\in\mathfrak{X}(M)$ be its tangential part along $\psi$, defined as follows: for every $q\in M$, $X_q$ is the orthogonal projection of $\Xbar_{\psi(q)}$ onto the tangent subspace $T_q M \subseteq T_{\psi(q)}\Mbar$.

\begin{lemma}
    Let $\psi:(M^m,g)\to(\Mbar^{m+1},\gbar)$ be an isometric immersion between semi-Riemannian manifolds, $\Xbar\in\mathfrak X(\Mbar)$ a vector field and $X\in\mathfrak X(M)$ its tangential part along $\psi$. For any couple of tangent vectors $V,W\in T_p M$, $p\in M$,
    \begin{equation} \label{covhyp}
        g(\nabla_V X,W) = \gbar(\nablabar_V \Xbar, W) + \gbar(\Xbar,\II(V,W))
    \end{equation}
    and for any choice of a local unit normal vector field $N$ on $M$,
    \begin{equation} \label{divhyp}
        \div(X) = \divbar(\Xbar) + \gbar(m\mathbf{H},\Xbar) - \gbar(N,N)\gbar(\nablabar_N \Xbar,N),
    \end{equation}
    where $\div$ and $\divbar$, respectively, are the divergence operators induced by $\nabla$ and $\nablabar$.
\end{lemma}

\begin{proof}
    Let $U,\Ubar,Z$ be as above. Then $X=\Xbar-\frac{\gbar(\Xbar,Z)}{\gbar(Z,Z)}Z$ on $M$, so
    \[
        g(\nabla_V X, W) = \gbar(\nablabar_V X, W) = \gbar(\nablabar_V \Xbar, W) - \gbar\left(\nablabar_V \left(\frac{\gbar(\Xbar,Z)}{\gbar(Z,Z)}Z\right), W\right).
    \]
    Since $\gbar(Z,W)=0$ and $\II(V,W)$ is a multiple of $Z$,
    \begin{align*}
        - \gbar\left(\nablabar_V \left(\frac{\gbar(\Xbar,Z)}{\gbar(Z,Z)}Z\right), W\right) & = - \frac{\gbar(\Xbar,Z)}{\gbar(Z,Z)} \gbar(\nablabar_V Z, W) = \frac{\gbar(\Xbar,Z)\gbar(\II(V,W),Z)}{\gbar(Z,Z)} \\
        & = \gbar(\Xbar,\II(V,W)).
    \end{align*}
    Now let us consider a $g$-orthonormal basis $\{e_i\}_{1\leq i\leq m}$ of $T_p M$. If $\gbar(Z,Z)=\pm1$, that is, if $N=Z$ is a unit normal vector for $\psi$ at $p$, then
    \begin{equation} \label{divXXbar}
        \begin{split}
            \div(X) & = \sum_{i=1}^m g(e_i,e_i) g(\nabla_{e_i}X,e_i), \\
            \divbar(\Xbar) & = \sum_{i=1}^m \gbar(e_i,e_i) \gbar(\nablabar_{e_i}\Xbar,e_i) + \gbar(N,N)\gbar(\nablabar_N \Xbar,N).
        \end{split}
    \end{equation}
    Then \eqref{divhyp} follows from \eqref{covhyp}, \eqref{divXXbar} and the definition \eqref{Hdef} of $\mathbf{H}$.
\end{proof}

In this note we are interested in the case where $(\Mbar,\gbar)$ is a Lorentzian manifold, $\Xbar$ is a timelike conformal vector field and $\psi:M\to\Mbar$ is a spacelike hypersurface, that is, $g$ is a Riemannian metric. Conformality of $\Xbar$ means that there exists $\eta\in C^{\infty}(\Mbar)$ such that $\mathcal L_{\Xbar}\gbar = 2\eta \gbar$ on $\Mbar$, that is, $\gbar(\nablabar_V \Xbar, W) = \eta\gbar(V,W)$ for every $V,W\in T_q\Mbar$, $q\in\Mbar$. Let us set the notation
\begin{equation}
    \alpha = \sqrt{-\gbar(\Xbar,\Xbar)}, \qquad \mathcal H = \frac{\eta}{\alpha}.
\end{equation}
If $\hat\psi : \Sigma^m \to \Mbar$ is a (necessarily spacelike) hypersurface such that $\Xbar$ is orthogonal to $\hat\psi$ at some point $x\in\Sigma$, then by \eqref{IIlie} the second fundamental form $\hat\II$ of $\hat\psi$ satisfies $\gbar(\hat\II(V,W),\Xbar) = -\eta\hat g(V,W)$ for every $V,W\in T_x\Sigma$, with $\hat g = \hat\psi^{\ast}\gbar$. $\hat N = \alpha^{-1}\hat X$ is a unit normal vector for $\hat\psi$ at $x$ and therefore $\hat\II(V,W) = \gbar(\hat N,\hat N)\gbar(\hat\II(V,W),\hat N)\hat N = \eta\alpha^{-1}\hat g(V,W)\hat N$. So, $\hat\II = \mathcal H \hat g \otimes \hat N$ at $x$, that is, $\hat\psi$ is umbilical at $x$ with mean curvature vector
\begin{equation}
    \hat{\mathbf{H}} = \mathcal H \hat N = \frac{\eta}{\alpha^2}\Xbar.
\end{equation}

The existence of a global timelike vector field $\Xbar$ implies that $\psi:M\to\Mbar$ is two-sided, that is, $M$ admits a global (timelike) unit normal vector field $N$. In fact, on every open subset $U\subseteq M$ admitting a local unit normal vector field $N_U : U \to TU^{\bot}$ the function $\gbar(N_U,\Xbar)$ is always nonzero since $N_U$ and $\Xbar$ are both timelike and therefore cannot be orthogonal at any point. So, consider a family $\{N_a\}_{a\in I}$ of local unit normal vector fields defined on the elements of an open cover $\{U_a\}_{a\in I}$ of $M$ and such that $\gbar(N_a,\Xbar)<0$ on $U_a$. The conditions $\gbar(N_a,N_a)=-1$ and $\gbar(N_a,\Xbar)<0$ uniquely determine $N_a$ at every point of $U_a$, so $N_a = N_b$ on $U_a \cap U_b$ for every $a,b\in I$ and therefore we can glue together these vectors to obtain a global unit normal vector field $N:M\to TM^{\bot}$ satisfying $\gbar(N,\Xbar)<0$ on $M$.

In the following, we will always assume that $N$ is chosen so that $\gbar(N,\Xbar)<0$, that is, with the same time-orientation of $\Xbar$. In this case, we will also say that $N$ is future-pointing. The mean curvature vector $\mathbf{H}$ then induces a mean curvature function $H = -\gbar(\mathbf H,N)\in C^{\infty}(M)$ for which $\mathbf H = HN$. By the wrong-way Cauchy-Schwarz inequality, $\gbar(N,\alpha^{-1}\Xbar)\leq -1$, so we can introduce the hyperbolic angle function $\theta\in C^{\infty}(M)$ via its hyperbolic cosine
\[
    \cosh\theta = - \gbar(N,\alpha^{-1}\Xbar) = - \frac{\gbar(N,\Xbar)}{\alpha}.
\]
In this setting, recalling \eqref{lieg} we can express formulas \eqref{covhyp} and \eqref{divhyp} as
\begin{equation} \label{covconf}
    g(\nabla_V X, W) = \eta g(V,W) + g(AV,W)\alpha\cosh\theta,
\end{equation}
with $A=-\nablabar_{(\,\cdot\,)}N$ the shape operator of $\psi$ induced by $N$, and
\begin{equation} \label{divconf}
    \div(X) = m\eta - mH\alpha\cosh\theta = m\alpha(\mathcal H - H\cosh\theta).
\end{equation}

\subsection{Doubly warped product spacetimes}

Examples of spacetimes admitting a timelike conformal vector field include doubly warped Lorentzian product spacetimes. As we have previously said, by a doubly warped product spacetime $\Mbar = \IxP$ we mean a product manifold $\Mbar = I \times \PP$, where $(\PP,\sigma)$ is a connected $m$-dimensional Riemannian manifold and $I\subseteq\mathbb R$ is an open interval, endowed with the Lorentzian metric $\gbar$ given at $(t,x)\in\Mbar$ by
\begin{equation}
	\label{metric}
	\gbar_{(t,x)} = -h(x)^2\pi_I^{\ast}(\di t^2) + \rho(t)^2\pi_{\PP}^{\ast}(\sigma)
\end{equation}
with $h\in C^{\infty}(\PP)$, $\rho\in C^{\infty}(I)$ positive functions and $\pi_I:\Mbar\to I$, $\pi_{\PP}:\Mbar\to \PP$ the projections onto the factors of $I\times P$. In the following, with an abuse of notation we write $\rho,\rho',h$ to denote the functions $\rho\circ\pi_I,\rho'\circ\pi_I,h\circ\pi_{\PP}\in C^{\infty}(\Mbar)$. A time orientation for $\Mbar$ is given by the timelike vector $\partial_t := \partial /\partial t$. Note that from these models we can reobtain a standard static spacetime by setting $\rho(t) \equiv 1$, as well as a Generalized Robertson-Walker spacetime when $h(x) \equiv 1$. 

\begin{lemma}
	\label{lemmaconfo}
	The timelike vector field $\Xbar = \rho\partial_t$ is conformal on $\Mbar = \IxP$ and
	\begin{equation} \label{IxPconf}
	    \mathcal L_{\Xbar}\gbar = 2\rho'\gbar.
	\end{equation}
\end{lemma}

\begin{proof}
	Let $(t_0,x_0)\in\Mbar$ be a given point and let $\{x^i\}$ be a coordinate system $(\PP,\sigma)$ on a neighbourhood $U\subseteq \PP$ of $x_0$. Then $\{t,x^i\}$ is a coordinate system for $\Mbar$ defined on $I\times U\ni(t_0,x_0)$ and the $(m+1)$-ple $\{e_1,\dots,e_{m+1}\} := \left\{ \rho^{-1}\partial_1,\dots,\rho^{-1}\partial_m,h^{-1}\partial_t\right\}$	is a local frame for $\Mbar$ on $I\times U$. A direct computation shows that $[\Xbar,e_{\mu}] = - \rho' e_{\mu}$ on $I\times U$ for $1 \leq \mu \leq m+1$. For every $1\leq \mu_1,\mu_2 \leq m+1$, the product $\gbar(e_{\mu_1},e_{\mu_2})$ is constant along the curve $I\times\{x_0\}$, so we have
	\begin{align*}
	(\mathcal L_{\Xbar}\gbar)(e_{\mu_1},e_{\mu_2}) & = \Xbar(\gbar(e_{\mu_1},e_{\mu_2})) - \gbar(e_{\mu_1},[\Xbar,e_{\mu_2}]) - \gbar(e_{\mu_2},[\Xbar,e_{\mu_1}]) \\
	& = 0 - \gbar(e_{\mu_1},-\rho' e_{\mu_2}) - \gbar(e_{\mu_2},-\rho' e_{\mu_1}) = 2\rho' \gbar(e_{\mu_1},e_{\mu_2})
	\end{align*}
	at $(t_0,x_0)$. Since $(\mathcal L_{\Xbar}\gbar)_{(t_0,x_0)}$ is a bilinear form on $T_{(t_0,x_0)}\Mbar$, \eqref{IxPconf} follows.
\end{proof}

A doubly warped product $\Mbar = \IxP$ is foliated by the level sets $\Sigma_t = \{t\}\times\PP$, $t\in I$, of the coordinate function $\pi_I : \Mbar \to I$. They are always orthogonal to $\Xbar$, so all of them are totally umbilical hypersurfaces by the previous discussion (see also Prop. 2.2 in \cite{U}). In this setting, we have
\[
    \alpha = \sqrt{-\gbar(\Xbar,\Xbar)} = h\rho, \qquad \mathcal H = \frac{\eta}{\alpha} = \frac{\rho'}{h\rho}.
\]
Now, let $\psi : M \to \IxP$ be a spacelike immersed hypersurface, that is, assume that $g=\psi^{\ast}\gbar$ is a Riemannian metric on $M$. Again, with a little abuse of notation, we denote with $\rho$, $\rho'$, $h$ the functions $\rho\circ\pi_I\circ\psi$, $\rho'\circ\pi_I\circ\psi$, $h\circ\pi_{\PP}\circ\psi\in C^{\infty}(M)$. Since the coordinate function $\pi_I$ has gradient $\nablabar\pi_I = -h^{-2}\partial_t = -\rho^{-1}h^{-2}\Xbar$, if we introduce the height function $\tau = \pi_I\circ\psi$ then the tangential part $X$ of $\Xbar$ along $\psi$ satisfies $- X = \rho h^2\nabla\tau$. Then, equation \eqref{divconf} can be restated as
\begin{equation}
	\label{divrt}
	\div(\rho h^2\nabla\tau) = mh\rho(H\cosh\theta - \mathcal H)
\end{equation}
or, equivalently,
\begin{equation} \label{Hvectorandflap}
	\Delta_{-\log(\rho h^2)}\tau = \frac{m}{h}(H\cosh\theta - \mathcal H)
\end{equation}
where, for $f\in C^1(M)$, $\Delta_f$ is the symmetric diffusion operator
\[
	\Delta_f = e^{f}\div(e^{-f}\nabla\quad) = \Delta - g(\nabla f,\nabla\quad).
\]
If we let $\mathcal R$ be an antiderivative of $\rho$ on $I$, \eqref{divrt} is also equivalent to
\begin{equation} \label{Hvectorandflap2}
	\Delta_{-\log(h^2)}(\mathcal R\circ\tau) = \frac{m\rho}{h}(H\cosh\theta - \mathcal H)
\end{equation}

\section{Rigidity of $H$-hypersurfaces with $H\geq\mathcal H$}

\label{se2}

The aim of this section is to prove some rigidity results for compact or complete spacelike hypersurfaces immersed in a spacetime $\Mbar$ carrying a timelike conformal vector field $\Xbar$. We shall refer again to the notation
\begin{equation} \label{alphaH}
    \mathcal L_{\Xbar}\gbar = 2\eta\gbar,  \qquad \alpha = \sqrt{-\gbar(\Xbar,\Xbar)}, \qquad \mathcal H = \frac{\eta}{\alpha}
\end{equation}
introduced in the previous section. We first consider the case of a doubly warped product spacetime $\Mbar = \IxP$, where $\Xbar=\rho\partial_t$ and
\begin{equation} \label{alphaHdw}
    \mathcal L_{\Xbar}\gbar = 2\rho'\gbar, \qquad \alpha = h\rho, \qquad \mathcal H = \frac{\rho'}{h\rho}.
\end{equation}

Let us remark some facts about spacelike immersions into doubly warped product spacetimes. Given $\psi: M \longrightarrow \overline{M}$ a spacelike hypersurface, we can define its projection on $\PP$ by $\varphi_{\PP} := \pi_{\PP} \circ \psi: M \longrightarrow \PP$. Note that
\begin{equation}
    \label{metcomp}
    \varphi_{\PP}^*(\sigma) = \rho^{-2}\psi^{\ast}(\gbar + h^2\di t^2) \geq \rho^{-2} \psi^{\ast}\gbar = \rho^{-2} g
\end{equation}
with $g=\psi^{\ast}\gbar$ the metric induced by $\psi$ on $M$. Since $\varphi_{\PP}$ is a local diffeomorphism, by Lemma 3.3 in Chapter 7 of \cite{DC}.
\begin{lemma}
	\label{lemcm}
	Let $\psi: M \longrightarrow \overline{M}$ be a complete spacelike hypersurface in a doubly warped product spacetime $\overline{M} = \IxP$ and let $\varphi_{\PP} := \pi_{\PP} \circ \psi$ be its projection on $\PP$. Then $\varphi_{\PP}$ is a local diffeomorphism which is a covering map when $\rho$ is bounded on $M$.
\end{lemma}

The following is a direct consequence of Lemma \ref{lemcm}.

\begin{proposition}
	\label{prpcvcp}
	Let $\Mbar = \IxP$ be a doubly warped product spacetime.
	\begin{enumerate}
		\item If $\overline{M}$ admits a compact spacelike hypersurface, then $\PP$ is compact.
		\item If the universal covering of $\PP$ is compact, then any complete spacelike hypersurface immersed in $\overline{M}$ where $\rho$ is bounded is compact.
	\end{enumerate}
\end{proposition}

We now use equation (\ref{divrt}) to deduce the following

\begin{theorem}
	\label{thcho}
	Let $\overline{M} = \IxP$ be a spatially closed doubly warped product spacetime, that is, assume that $\PP$ is compact. Then, $\overline{M}$ admits a compact maximal hypersurface if and only if it admits a totally geodesic spacelike slice.
\end{theorem}

\begin{proof}
	Any totally geodesic spacelike slice of a spatially closed doubly warped product is clearly a maximal compact hypersurface. Vice versa, let $\psi:M\to\IxP$ be a compact hypersurface with mean curvature $H$. From (\ref{divrt}) and \eqref{alphaHdw} we get
	\begin{equation}
	\label{divrt2}
	    \div\left(\rho h^2\nabla\tau\right) =  m \rho h H \cosh \theta - m \rho'.
	\end{equation}
	Integrating (\ref{divrt2}) on a compact hypersurface $M$ we infer the integral formula
	\begin{equation}
	\label{inte}
	    \int_M \left( \rho h H \cosh \theta - \rho' \right) = \int_M \div(\rho h^2\nabla\tau) = 0.
	\end{equation}
	If $M$ is maximal, then $H\equiv0$ and (\ref{inte}) reduces to
	$$\int_M \rho' =0.$$
	Thus, there must exist $t_0 \in I$ such that $\rho'(t_0) = 0$. Since the spacelike slices are totally umbilic and $\mathcal H(t_0,\;\cdot\;)\equiv0$, $\Sigma_{t_0} = \{t_0\}\times\PP$ is a totally geodesic spacelike slice.
\end{proof}

The above theorem extends to the present ambient spacetimes a result obtained by Choquet-Bruhat in \cite{CB} for Robertson-Walker spacetimes (see also Prop. 4.1 in \cite{ARS}). As a second consequence of equation \eqref{divrt} we prove the next rigidity result, which generalizes Theorem 1 of \cite{PRR} to doubly warped product spacetimes.

\begin{theorem}
	\label{thcpt}
	Let $\Mbar = \IxP$ be doubly warped spacetime satisfying $\rho'\geq0$ on $I$. If $\psi:M\to\Mbar$ is a connected compact spacelike hypersurface whose mean curvature in the direction of the future-pointing normal satisfies
	\[
	    H \geq \mathcal{H}\circ\psi = \frac{\rho'}{h\rho}\circ\psi
	\]
	then $\psi(M)$ is a spacelike slice.
\end{theorem}

\begin{proof}
	Let $\psi : M \to \Mbar$ be a compact spacelike hypersurface satisfying $H\geq\mathcal H$. By equation \eqref{divrt}, we have
	\begin{equation} \label{div<=0}
	    \div(\rho h^2\nabla\tau) = mh\rho(H\cosh\theta - \mathcal H) \geq m\rho\mathcal H(\cosh\theta - 1) \geq 0
	\end{equation}
	since $\mathcal H\geq0$. Integrating \eqref{div<=0} on the compact manifold $M$ and applying the divergence theorem we get
	\[
        0 = \int_M \div(\rho h^2\nabla\tau)
    \]
	so $\div(\rho h^2\nabla\tau)=0$ on $M$. Then $\tau$ is constant on $M$ by the strong maximum principle. This is clearly equivalent to saying that $\psi(M)$ is contained in a slice $\Sigma_{t_0} = \{t_0\} \times \PP$ for some $t_0\in I$. Since $M$ is compact, $\psi(M)\subseteq\Sigma_{t_0}$ is closed. The map $\psi : M \to \Sigma_{t_0}$ is an immersion between manifolds of equal dimension, so it is a local diffeomorphism and therefore an open map. Hence, the nonempty image $\psi(M)$ is both open and closed in the connected slice $\Sigma_{t_0}$ and we conclude that $\psi(M) = \Sigma_{t_0}$.
\end{proof}

Clearly, the conclusion of Theorem \ref{thcpt} is still true if we replace the assumption $H\geq\mathcal H\circ\psi\geq0$ with $H\leq\mathcal H\circ\psi\leq0$. A similar property holds in the more general case of a spacetime admitting a timelike conformal vector field, provided that the stronger condition $H\geq\mathcal H\circ\psi>0$ is satisfied.

\begin{theorem}
    \label{thcptgen}
    Let $(\Mbar^{m+1},\gbar)$ be a spacetime carrying a conformal timelike vector field $\Xbar\in\mathfrak{X}(M)$ and suppose that $\mathcal L_{\Xbar}\gbar = 2\eta\gbar$ with $\eta>0$. Let $\psi:M^m\to\Mbar$ be a compact immersed spacelike hypersurface whose mean curvature function in the direction of the future-pointing normal satisfies $H \geq \mathcal H\circ\psi$. Then $\Xbar$ is orthogonal to $\psi$ and $H=\mathcal H\circ\psi$ on $M$.
\end{theorem}

\begin{proof}
    By equation \eqref{divconf}, the tangential part $X$ of $\Xbar$ along $\psi$ satisfies
    \begin{equation} \label{div<0}
        \div(X) = -m\alpha(H\cosh\theta-\mathcal H) \leq -m\alpha\mathcal H(\cosh\theta-1) \leq 0,
    \end{equation}
    because $\mathcal H = \alpha^{-1}\eta > 0$. As in the proof of Theorem \ref{thcpt}, an application of the divergence theorem yields $\div(X)\equiv0$ on $M$. Since $\mathcal H > 0$, by \eqref{div<0} we conclude that $\cosh\theta\equiv1$ on $M$.
\end{proof}

\begin{remark}
\normalfont
    Note that in the hypotheses of Theorem \ref{thcptgen}, the orthogonal distribution of $\Xbar$ is not assumed to be integrable, that is, the \textit{a priori} existence of immersed hypersurfaces orthogonal to $\Xbar$ is not assumed. Also note that the conclusion of the theorem is false, in general, if vanishing of $\eta$ is allowed. For instance, the Lorentzian surface $(\Mbar,\gbar)$ obtained by endowing the cylinder $\Mbar = \{(x,y,t)\in\R^3 : x^2 + y^2 = 1\}$ with the Lorentzian metric $\gbar = -\di t^2 + \di x^2 + \di y^2$ induced by the restriction to $\Mbar$ of the Lorentz-Minkowski metric of $\R^3$ is foliated by compact spacelike geodesics $\psi_{t_0} : \mathbb{S}^1 \to \Mbar : \theta \mapsto (\cos\theta,\sin\theta, t_0)$, $t_0\in\R$. However, $\Mbar$ carries a family of timelike vector fields
    \[
        \Xbar_a = \partial_t - ay\partial_x + ax\partial_y, \qquad -1 < a < 1
    \]
    satisfying $\mathcal L_{\Xbar_a}\gbar = 0$ for every $-1 < a < 1$, but, for every $t_0\in\R$, $\psi_{t_0}$ is orthogonal to $\Xbar_a$ if and only if $a=0$. In fact, when $a\neq0$, the maximal codimension $1$ spacelike submanifolds orthogonal to $\Xbar_a$ are the noncompact geodesics $\gamma_{a,t_0} : \R \to \Mbar : s \mapsto (\cos s, \sin s, as + t_0)$, $t_0\in\R$.
\end{remark}

Note that when $\Mbar$ is a doubly warped product spacetime, the values of $H$ on a compact spacelike hypersurface naturally relate with those of $\mathcal H$ as a consequence of the maximum principle. In fact, we have the following

\begin{theorem} \label{thmhbound}
	Let $\psi: M \longrightarrow \overline{M}$ be a spacelike hypersurface in a doubly warped product spacetime $\Mbar = \IxP$. Suppose that there exist two points $p_0$ and $p^0$, respectively, where the height function $\tau$ attains local minimum and maximum values. Then,
	\begin{equation} \label{hbo}
		\mathcal{H}(\psi(p_0)) \leq H(p_0), \qquad \mathcal{H}(\psi(p^0)) \geq H(p^0).
	\end{equation}
\end{theorem}

\begin{proof}
	Since $p_0$ and $p^0$ are locally extremal for $\tau$, at these 
	points we have $\nabla \tau = 0$, so $\rho h^2 \Delta \tau = \diver(\rho h^2 \nabla \tau)$ and $X=0$, implying that $N = \alpha^{-1}\Xbar$ and therefore $\cosh\theta=1$. Since $p_0$ is a local minimum point for $\tau$, we have $\Delta\tau(p_0) \geq 0$ and then
	\[
	    0 \leq \frac{\div(\rho h^2\nabla\tau)}{mh\rho} = H - \mathcal{H} \qquad \text{at } \, p_0.
	\]
	Similarly, we deduce $0 \geq H - \mathcal{H}$ at $p^0$.
\end{proof}

When $H$ is constant on $M$ and $\rho\equiv1$, that is, when $\Mbar$ is a standard static spacetime, the inequalities \eqref{hbo} are clearly satisfied if and only if $H\equiv0$.

\begin{corollary}
	\label{coross}
	Let $\Mbar = I {}_h\hspace{-2pt}\times \PP$ be a standard static spacetime, $\psi:M\to\Mbar$ a connected spacelike hypersurface with constant mean curvature. If the height function attains both locally minimal and maximal values on $M$, then $\psi(M)$ is contained in a totally geodesic spacelike slice.
\end{corollary}

\begin{proof}
    Since $H$ is constant on $M$ and $\mathcal H \equiv 0$ on $\Mbar$, by Theorem \ref{thmhbound} it must be $H\equiv0$ on $M$. So $\tau$ satisfies $\Delta_{-\log(h^2)}\tau = 0$ on $M$ and attains a local maximum at some point of $M$. By the strong maximum principle together with connectedness of $M$ and the unique continuation property for the equation $\Delta_{-\log(h^2)}u = 0$, $\tau$ is constant on $M$. So, $\psi(M)$ is contained in a spacelike slice.
\end{proof}

We conclude by giving different versions of Theorem \ref{thcpt} in the complete noncompact case, replacing compactness of $M$ with different assumptions. The first one relies on the following result due to Al\'ias, Caminha, do Nascimento, see Theorem 2.1 of \cite{ACdoN}.

\begin{proposition} \label{divinf}
    Let $(M,g)$ be a complete, noncompact Riemannian manifold, $V\in\mathfrak X(M)$. Assume that there exists $f\in C^{\infty}(M)$, $f\geq0$, $f\not\equiv0$ such that
    \begin{equation}
        g(X,\nabla f) \geq 0 \quad \text{on } \, M
    \end{equation}
    and
    \begin{equation} \label{phito0}
        f(x) \to 0 \quad \text{as } \, x \to \infty \, \text{ in } \, M.
    \end{equation}
    If $\div(X)\geq0$ then
    \begin{equation} \label{Xvan}
        i) \; \; \div(X) \equiv 0 \quad \text{on } \, M \setminus f^{-1}(0) \qquad \text{and} \qquad ii) \; \; g(X,\nabla f) \equiv 0 \quad \text{on } \, M.
    \end{equation}
\end{proposition}

\begin{remark}
\normalfont
    If $\div$ is replaced everywhere in Proposition \ref{divinf} by the weighted divergence operator $\div_{\varphi}$ defined by $\div_{\varphi}(X) = \div(X) - g(\nabla\varphi,X) = e^{\varphi} \div(e^{-\varphi} X)$, with $\varphi\in C^{\infty}(M)$, then the resulting statement is also true. In fact, it suffices to apply the above Proposition to the vector field $e^{-\varphi}X$.
\end{remark}

\begin{theorem} \label{thasymp}
    Let $\psi:M\to\Mbar=\IxP$ be a spacelike complete, noncompact, connected hypersurface in a doubly warped product spacetime satisfying $\rho'\geq0$. Let $t_0\in I$ and suppose that $\psi(M)$ is above the slice $\Sigma_{t_0}$ and asymptotic to it at infinity. If $H \geq \mathcal H\circ\psi$ on $M$, then $\psi(M)=\Sigma_{t_0}$.
\end{theorem}

\begin{proof}
    As in the proof of Theorem \ref{thcpt}, the assumption $H\geq\mathcal H$ together with equation \eqref{divrt} yields
    \[
        \div(\rho h^2\nabla\tau) = mh\rho(H\cosh\theta - \mathcal H) \geq 0
    \]
    for $\tau=\pi_I\circ\psi$ the height function of $\psi$. Since $\psi(M)$ is above the slice $\Sigma_{t_0}=\{t_0\}\times\PP$ we have $\tau \geq t_0$ on $M$. Furthermore, $\tau(x)\to t_0$ as $x\to\infty$ in $M$. Thus the function $f=\tau-t_0$ satisfies $f\geq0$ on $M$ and, for $X=\rho h^2\nabla\tau$,
    \[
        g(X,\nabla f) = \rho h^2|\nabla\tau|^2 \geq 0.
    \]
    We reason by contradiction and we suppose that $\psi(M)\neq\Sigma_{t_0}$. Then for some $x\in M$ we have $\varphi(x)>0$, so that $\varphi\not\equiv0$ on $M$. We can thus apply Proposition \ref{divinf} to deduce that
    \[
        0 \equiv g(X,\nabla f) = \rho h^2|\nabla\tau|^2 \quad \text{on } \, M.
    \]
    Since $\rho h^2>0$ on $M$, it follows that $\varphi\equiv0$ on $M$, contradiction.
\end{proof}

In the next two results, the role of compactness is played by the parabolicity of a certain operator.

\begin{theorem} \label{thparab1}
	Let $\Mbar = \IxP$ be a doubly warped product spacetime satisfying $\rho'\geq0$ on $I$ and let $\psi : M \to \Mbar$ be a spacelike complete, noncompact hypersurface. Suppose that
	\begin{equation} \label{hrhonoL1}
		\int_1^{+\infty} \left( \int_{\partial B_r(o)} \rho h^2 \right)^{-1} \di r = +\infty
	\end{equation}
	for some reference point $o\in M$. If $H \geq \mathcal H\circ\psi$ on $M$ and the height function $\tau$ is bounded above, then $\psi(M)$ is a spacelike slice.
\end{theorem}

\begin{proof}
	Under the assumption $H\geq\mathcal H\geq0$, equation \eqref{Hvectorandflap} yields
	\begin{equation} \label{Hcomp1}
		\Delta_{-\log \rho h^2} \tau = \frac{m}{h}(H\cosh\theta - \mathcal H) \geq \frac{m\mathcal H}{h}(\cosh\theta - 1) \geq 0.
	\end{equation}
	Next completeness of $M$ and \eqref{hrhonoL1} imply that the operator $\Delta_{-\log \rho h^2}$ is parabolic on $M$, see Chapter 4 in \cite{AMR}. Thus since $\tau$ is bounded above from \eqref{Hcomp1} we deduce that $\tau$ is constant, and this implies that $\psi(M)$ is contained in a spacelike slice. By completeness of $M$, $\psi(M)$ must in fact be a slice.
\end{proof}

Observe that condition \eqref{hrhonoL1} involves both $\rho$ and $h$. We can get rid of $\rho$ by considering equation \eqref{Hvectorandflap2} instead of \eqref{Hvectorandflap}. If $\mathcal R$ is an antiderivative of $\rho$ on $I$, then $\mathcal R' = \rho > 0$ and therefore
\[
	\sup_{M} \mathcal R(\tau) = \lim_{t\to\sup_M \tau} \mathcal R(t)
\]
when $\tau$ is the height function of an immersion $\psi:M\to\Mbar$. So, for instance, if $\sup_M\tau \in I$ then $\mathcal R(\tau)$ is bounded above. By applying the argument used in the proof of Theorem \ref{thparab1} to equation \eqref{Hvectorandflap2} we have the following Theorem. Note that condition \eqref{hnoL1} below is equivalent to \eqref{hrhonoL1} if $\rho$ is bounded above and stays away from $0$, but otherwise the two seem independent.

\begin{theorem} \label{thparab2}
	Let $\Mbar = \IxP$ be a doubly warped product spacetime satisfying $\rho'\geq0$ on $I$. Let $\psi : M \to \Mbar$ be a spacelike complete hypersurface such that
	\begin{equation} \label{hnoL1}
		\int_1^{+\infty} \left( \int_{\partial B_r(o)} h^2 \right)^{-1} \di r = +\infty
	\end{equation}
	for some reference point $o\in M$. If $H \geq \mathcal H\circ\psi$ on $M$ and $\sup_M \tau \in I$, then $\psi(M)$ is a spacelike slice.
\end{theorem}

\section*{Acknowledgements}

The second author is partially supported by Spanish MINECO and ERDF project MTM2016-78807-C2-1-P.

%
%

\begin{thebibliography}{99.}%
%
%

\bibitem{ACdoN} Al\'ias, L.J., Caminha, A., do Nascimento, Y.: A maximum principle at infinity with applications to geometric vector fields. J. Math. Anal. Appl. (in press) doi: 10.1016/j.jmaa.2019.01.042

\bibitem{AMR} Al\'\i as, L.J., Mastrolia, P., Rigoli, M.: Maximum principles and geometric applications. Springer (2016)

\bibitem{ARS} Al\'\i as, L.J., Romero, A., S\'anchez, M.: Uniqueness of complete spacelike hypersurfaces of constant mean curvature in Generalized Robertson-Walker spacetimes. Gen. Relat. Gravit. \textbf{27}, 71--84 (1995)

\bibitem{Al} Allison, D.E.: Lorentzian warped products and static space-times. PhD Thesis, University of Missouri-Columbia (1985)

\bibitem{BP} Beem, J.K, Powell, T.G.: Geodesic completeness and maximality in Lorentzian warped products. Tensor (N.S.) \textbf{39}, 31--36 (1982)

\bibitem{BS} Bernal, A.N., S\'anchez, M.: On smooth Cauchy hypersurfaces and Geroch's splitting theorem. Commun. Math. Phys. \textbf{257}, 43--50 (2005)

\bibitem{BO} Bishop, R.L., O'Neill, B.: Manifolds of negative curvature. Trans. Amer. Math. Soc. \textbf{145}, 1--49 (1969)

\bibitem{CRR} Caballero, M., Romero, A., Rubio, R.M.: Constant mean curvature spacelike hypersurfaces in Lorentzian manifolds with a timelike gradient conformal vector field. Classical Quant. Grav. \textbf{28}:145009, 13 pp. (2011)

\bibitem{CPR} Colombo, G., Pelegr\'in, J.A.S., Rigoli, M.: Spacelike hypersurfaces 
in standard static spacetimes. Gen. Relat. Gravit. \textbf{51}:1, 66 pp. (2019)

\bibitem{CPR2} Colombo, G., Pelegr\'in, J.A.S., Rigoli, M.: Stable maximal hypersurfaces in Lorentzian spacetimes. Nonlinear Anal. \textbf{179}, 354--382 (2019)

\bibitem{CB} Choquet-Bruhat, Y.: Maximal submanifolds and submanifolds with constant mean extrinsic curvature of a lorentzian manifold. Ann. Scuola Norm.-Sci. \textbf{3}, 361--376 (1976)

\bibitem{DC} Do Carmo, M.P.: Riemannian Geometry. Birkh\"auser (1992)

\bibitem{Ea} Eardley, D., Isenberg, J., Marsden, J., Moncrief, V.: Homothetic and conformal symmetries of solutions to Einstein's equations. Comm. Math. Phys. \textbf{106}, 137--158 (1996)

\bibitem{Ger} Geroch, R.: Domain of dependence. J. Math. Phys. \textbf{11}, 437--449 (1970)

\bibitem{GO} Guti\'errez, M., Olea, B.: Semi-Riemannian manifolds with a doubly warped structure. Rev. Mat. Iberoam. \textbf{28}, 1--24 (2011)

\bibitem{O'N} O'Neill, B.: Semi-Riemannian Geometry with applications to Relativity. Academic Press, New York (1983)

\bibitem{PRR} Pelegr\'in, J.A.S., Romero, A., Rubio, R.M.: On uniqueness of the foliation by comoving observers restspaces of a Generalized Robertson-Walker spacetime. Gen. Relat. Gravit. \textbf{49}:16, 14 pp. (2017)

\bibitem{Ri} Ringstr\"om, H.: The Cauchy problem in General Relativity. ESI Lectures in Mathematics and Physics, European Mathematical Society (2009)

\bibitem{U} \"Unal, B.: Doubly warped product. Differ. Geom. Appl. \textbf{15}, 253--263 (2001)

\end{thebibliography}
%

\end{document}